\documentstyle[amsfonts,11pt]{article}

\newtheorem{theo}{Theorem}[section]
\newtheorem{defin}[theo]{Definition}
\newtheorem{remar}[theo]{Remark}

\newtheorem{lemma}[theo]{Lemma}

\newtheorem{Example}[theo]{Example}

\newcommand{\fdim}{\hspace*{\fill}$\Box$}
\newcommand{\dimostr}{{\bf Proof: }}
\newcommand{\real}{\Bbb{R}}
\newcommand{\R}{\Bbb{R}}

\newcommand{\complex}{\Bbb{C}}
\newcommand{\C}{\Bbb{C}}

\newcommand{\N}{\Bbb{N}}

\newcommand{\K }{K\"{a}hler}

\begin{document}

\centerline {\LARGE\bf Symplectic maps
of complex domains}
\vspace{0.3cm}
\centerline {\LARGE\bf  into complex space forms
 \footnote{During the preparation of this article the authors were supported
by the M.I.U.R. Project \lq\lq Geometric Properties of Real and
Complex Manifolds''. }}

\vspace{0.5cm}

\centerline{\small Andrea Loi} \centerline{\small Dipartimento di
Matematica e Informatica -- Universit\`{a} di Cagliari -- Italy}
\centerline{\small e-mail address: loi@unica.it}

\vspace{0.3cm} \centerline{\small and} \vspace{0.3cm}

\centerline{\small Fabio Zuddas} \centerline{\small Dipartimento
di Matematica e Informatica -- Universit\`{a} di Cagliari --
Italy} \centerline{\small e-mail address: fzuddas@unica.it}

\vspace{0.3cm}

\begin{abstract}
Let $M\subset{\complex}^n$ be a complex domain of ${\complex}^n$
endowed with a rotation invariant \K\ form $\omega_{\Phi}=
\frac{i}{2} \partial\bar\partial\Phi$. In this paper we describe
sufficient conditions on the \K\ potential $\Phi$ for $(M,
\omega_{\Phi})$ to admit a symplectic embedding (explicitely described in terms of $\Phi$)
into a complex space form of the same dimension of $M$. In particular we also provide
conditions on $\Phi$ for $(M, \omega_{\Phi})$ to admit global symplectic
coordinates. As an application of our results we  prove that each of the Ricci flat
(but not flat) \K\ forms on ${\complex}^2$ constructed by LeBrun
in \cite{lebrun1} admits explicitely computable global symplectic
coordinates.


{\it{Subj.Class}}: 53C55, 58C25, 53D05, 58F06.

\end{abstract}

\section{Introduction and statements of the main results}
Let $(M, \omega)$  and $(S, \Omega)$ be two symplectic manifolds
of dimension  $2n$ and $2N$, $n\leq N$,  respectively. Then, one
has the following natural and fundamental question.

\vskip 0.1cm

\noindent
{\bf Question 1.}
{\em Under which conditions there exists a
symplectic embedding  $\Psi: (M, \omega )\rightarrow (S, \Omega )$,
namely a smooth embedding $\Psi: M\rightarrow S$
satisfying $\Psi^*(\Omega)=\omega$?}

\vskip 0.1cm

Theorems A, B, C and D below give a topological answer to the
previous question when   $ \Omega$ is the \K\  form  of  a
$N$-dimensional complex space form $S$, namely $(S, \Omega)$
is either the complex Euclidean space $({\complex}^N, \omega_0)$,
the complex hyperbolic space $({\complex}H^N, \omega_{hyp})$
or the complex projective space $({\complex}P^N, \omega_{FS})$
(see below for the definition of the symplectic (\K ) forms  $\omega_0$,
$\omega_{hyp}$ and $\omega_{FS}$ ). Indeed these theorems are
consequences of  Gromov's h-principle   \cite{gr2}
(see  also Chapter 12 in \cite{el} for a beautiful Êdescription of
Gromov's work ).

\vskip 0.3cm

\noindent
{\bf Theorem A (Gromov \cite{gr2}, see also \cite{gr0})}
{\em Let $(M, \omega)$ be a contractible  symplectic manifold.
Then there exist a non-negative integer $N$
 and a symplectic embedding
 $\Psi :(M, \omega)\rightarrow ({\complex}^{N}, \omega_0)$,
where
$\omega_0=\sum_{j=1}^Ndx_j\wedge dy_j$
denotes the standard symplectic form on ${\complex}^N={\real}^{2N}$. }

\vskip 0.3cm

 This was further  generalized by Popov as follows.

 \vskip 0.3cm

\noindent
{\bf Theorem B  (Popov \cite{popov})}
{\em Let $(M, \omega)$ be a   symplectic manifold.
Assume $\omega$ is exact,
namely $\omega =d\alpha$, for
 a $1$-form $\alpha$. Then there exist a non-negative integer $N$
 and a symplectic embedding
 $\Psi :(M, \omega)\rightarrow ({\complex}^{N}, \omega_0)$. }

  \vskip 0.3cm

Observe that the complex hyperbolic space
$({\complex}H^N,\omega_{hyp})$,
namely the unit ball
${\complex}H^N=\{z=(z_1,\dots, z_N)\in {\complex}^N|\ \sum_{j=1}^N|z_j|^2< 1\}$
in ${\complex}^N$ endowed with the hyperbolic form
 $ \omega_{hyp}=-\frac{i}{2}\partial\bar\partial\log (1-\sum_{j=1}^N|z_j|^2)$
is globally symplectomorphic to $({\complex}^N, \omega_0)$ (see
(\ref{s1}) in Lemma \ref{lemmadual} below) hence Theorem $B$
immediately implies

 \vskip 0.3cm

\noindent
{\bf Theorem C}
{\em Let $(M, \omega)$ be a   symplectic manifold.
Assume $\omega$ is exact.
 Then there exist a non-negative integer $N$
 and a symplectic embedding
 $\Psi :(M, \omega)\rightarrow ({\complex}H^{N}, \omega_{hyp})$. }

  \vskip 0.3cm

The following theorem, further generalized by Popov \cite{popov}
 to the noncompact case,  deals with the  complex projective ${\complex}P^{N}$,
equipped with the Fubini--Study form $\omega_{FS}$. Recall that if
$Z_0,\dots ,Z_N$ denote the  homogeneous
 coordinates on ${\complex}P^N$,
then, in   the affine chart $Z_0\neq 0$
endowed with coordinates  $z_j= \frac{Z_j}{Z_0}, j=1,\dots , N$,
the Fubini-Study form reads as
$$\omega_{FS}=\frac{i}{2}\partial\bar\partial\log  (1+\sum_{j=1}^N|z_j|^2).$$

 \newpage

\noindent {\bf Theorem D (Gromov \cite{gr0}, see also Tischler
\cite{ti}}) {\em Let $(M, \omega)$ be a compact symplectic
manifold such that $\omega$ is integral. Then there exist a
non-negative integer $N$
 and a symplectic embedding
 $\Psi :(M, \omega)\rightarrow ({\complex}P^{N}, \omega_{FS})$.}

 \vskip 0.3cm
 At this point a natural problem  is that to find the smallest
 dimension  of the complex space form where a given symplectic manifold $(M, \omega)$ can
 be symplectically embedded.
 In particular one can study the case of equidimensional symplectic maps,
 as expressed by the following interesting  question.

 \vskip 0.2cm

\noindent
{\bf Question 2.}
{\em Given a $2n$-dimensional symplectic manifold
$(M, \omega)$  under which conditions there exists a symplectic embedding $\Psi$
of  $(M, \omega)$ into $({\complex}^n, \omega_0)$ or $({\complex}P^n, \omega_{FS})$?}

 \vskip 0.2cm

\noindent
 Notice that locally there are not obstructions to the existence of such
 $\Psi$.
 Indeed, by a well-known theorem of Darboux
for every point $p\in M$ there exist a neighbourhood $U$ of $p$
and an embedding $\Psi:U\rightarrow {\real}^{2n} = {\complex}^n$
such that $\Psi^*(\omega_0)=\omega$.
In order to get a local embedding into
$({\complex}P^n, \omega_{FS})$  we can
assume (by shrinking $U$ if necesssary)
that $\Psi(U) \subset {\complex}H^n$.
 Therefore $f \circ
\Psi: U \rightarrow ({\complex}^n, \omega_{FS}) \subset
({\complex}P^n, \omega_{FS})$, with $f$ given by Lemma
\ref{lemmadual} below, is the desired embedding
satisfying $(f \circ\Psi)^*(\omega_{FS})=\Psi^*(\omega_0)=\omega$.
Observe also that Darboux's theorem is a special case of the following

\vskip 0.3cm

\noindent {\bf Theorem E (Gromov \cite{grihes})} {\em A
$2n$-dimensional symplectic manifold $(M, \omega)$ admits a
symplectic immersion into $({\complex}^n, \omega_0)$ if and only if the
following three conditions are satisfied: a) $M$ is open, b) the
form $\omega$ is exact, c) the tangent bundle $(TM, \omega )$ is a
trivial $Sp(2n)$-bundle. (Observe that a), b), c) are satisfied if
$M$ is contractible).}

 \vskip 0.3cm

It is worth pointing out that the previous theorem is not of any
help in order to attack Question 2 due to the existence of exotic
symplectic structures on ${\real}^{2n}$ (cfr. \cite{gr1}). (We
refer the reader to  \cite{ba} for an explicit construction of a
$4$-dimensional symplectic manifold diffeomorphic to ${\real}^4$
which cannot be symplectically embedded in $({\real}^4,
\omega_0)$).

\vskip 0.1cm In the case when our symplectic manifold $(M,
\omega)$
 is a \K\ manifold, with associated \K\ metric  $g$,
one can try to impose Riemannian  or holomorphic conditions  to
answer  the previous question. From the Riemannian point of view
the only complete and known result (to the authors' knowledge) is
the following global version of   Darboux's theorem.

 \vskip 0.3cm

\noindent
{\bf Theorem F (McDuff \cite{mc})}
{\em Let  $(M , g)$ be a simply-connected and complete
$n$-dimensional \K\ manifold of non-positive sectional curvature. Then there exists
a diffeomorphism $\Psi :M\rightarrow {\real}^{2n}$ such that $\Psi^*(\omega_{0})=\omega$.}

 \vskip 0.3cm

\noindent (See also \cite{cr1}, \cite{cr2},  \cite{cr3} and
\cite{dual} for further properties of McDuff's symplectomorphism).

\smallskip

The aim of this paper is to give  an answer to Question 2 in terms
of the \K\ potential of  the \K\ metric of  complex domains  (open
and connected) $M\subset {\complex}^n$ equipped with a
\K\ form $\omega$ which admits a  rotation
invariant \K\ potential.
 More precisely, throughout this paper
we  assume that there exists a  \K\ potential for
$\omega$,  namely a smooth function $\Phi:M\rightarrow
{\real}$ such that
$\omega=\frac{i}{2}\partial\bar\partial\Phi$,   depending  only on $|z_1|^2, \dots , |z_n|^2$,
where $z_1, \dots , z_n$ are the standard complex  coordinates on ${\complex}^n$.
Therefore, there exists a smooth function
$\tilde\Phi:\tilde M\rightarrow
{\real}$, defined on the open subset   $\tilde M\subset {\real}^n$ given by
\begin{equation}\label{tildem}
\tilde M=\{x=(x_1, \dots ,x_n)\in {\real}^n|\ x_j=|z_j|^2, z=(z_1, \dots z_n)\in M \}
\end{equation}
such
that
$$\Phi(z_1, \dots  ,z_n)=\tilde\Phi(x_1,\dots ,x_n), \ x_j=|z_j|^2,\ j=1,\dots ,n.$$
We  set $\omega :=\omega_{\Phi}$ and call
$\omega_{\Phi}$  a {\em rotation invariant} symplectic (\K ) form
with {\em associated} function $\tilde\Phi$. It is worth pointing
out that many interesting examples of \K\ forms on complex domains
are rotation invariant (even radial, namely depending only on
$r=|z_1|^2+\cdots +|z_n|^2$), since they often arise from
solutions of ordinary differential equations on the variable $r$
(see Example \ref{esempioloc2} below and also \cite{ca1} in the case of
extremal metrics).

\vskip 0.3cm

Our first  result is Theorem \ref{cor3} below where we  describe
explicit conditions in terms of the potential $\Phi$ for the
existence  of  an  explicit symplectic embedding of a rotation
invariant domain
 $(M, \omega_{\Phi})$ into a given complex space form $(S, \omega_{\Xi})$
 of the same dimension.
 In particular we find conditions on $\Phi$ for the existence of global symplectic
 coordinates of $(M, \omega_{\Phi})$.

\begin{theo}\label{cor3}
Let $M \subseteq {\C}^n$ be a complex domain
such that condition
\begin{equation}\label{intersezassi}
M \cap \{ z_j = 0 \} \neq \emptyset, \, j=1, \dots ,n
\end{equation}
is satisfied \footnote{Obviously  (\ref{intersezassi}) is
satisfied if $0\in M$, but there are other interesting cases, see
Examples
 \ref{esempioloc1} and \ref{esempioloc2} below,
 where this condition is fulfilled.}
 and  let $\omega_\Phi = \frac{i}{2}
\partial \bar{\partial} \Phi$ be a rotation invariant \K\ form
on $M$ with associated function $\tilde\Phi: {\tilde M} \rightarrow
{\R}$. Then

\begin{enumerate}
\item[(i)] there exists a uniquely  determined special
 \footnote{See (\ref{rotinv}) in  the next section
 for the definition of special map between complex domains.}
 symplectic immersion
 $$\Psi_0: (M, \omega_{\Phi}) \rightarrow ({\C}^n, \omega_0)$$
(resp. $\Psi_{hyp}: (M, \omega_{\Phi}) \rightarrow ({\C}H^n, \omega_{hyp})$) if and only if,
\begin{equation}\label{cond0}
\frac{\partial \tilde\Phi}{\partial x_k} \geq 0, \ \  k=1, \dots
,n.
\end{equation}
\item[(ii)]
there exists a uniquely determined special symplectic immersion
$$\Psi_{FS}:  (M, \omega_{\Phi}) \rightarrow
({\C}^n, \omega_{FS}),$$ if and only if
\begin{equation}\label{conda}
\frac{\partial \tilde\Phi}{\partial x_k} \geq 0, \ \  k=1, \dots
,n \ {\mbox and} \ \   \sum_{j=1}^n \frac{\partial
\tilde\Phi}{\partial x_j} x_j < 1,
\end{equation}
where we are looking at ${\C}^n \stackrel{i}{\hookrightarrow}
{\complex}P^n$ as the affine chart $Z_0\neq 0$ in ${\complex}P^n$
endowed with the restriction of the Fubini--Study form $\omega_{FS}$.
\end{enumerate}
Moreover, assume that  $0 \in M$.
If (\ref{cond0}) (resp.(\ref{conda}))  is satisfied
 then
 $\Psi_0$ (resp. $\Psi_{FS}$) is a
 global symplectomorphism
  (and hence $i \circ \Psi_{FS}:M \rightarrow {\complex}P^n$ is a
  symplectic embedding) if and only if
\begin{equation}\label{genconda}
\frac{\partial \tilde\Phi}{\partial x_k} > 0,  \ \  k=1, \dots ,n
\end{equation}
 and
 \begin{equation}\label{gencondb}
 \lim_{x \rightarrow \partial M} \sum_{j=1}^n
\frac{\partial \tilde\Phi}{\partial x_j} x_j = + \infty\  \
(resp.  \lim_{x \rightarrow \partial M} \sum_{j=1}^n
\frac{\partial \tilde\Phi}{\partial x_j} x_j = 1).\  \
 \footnote{For
 a rotation invariant continuous map $F: M \rightarrow {\R}$
we write
 $$\lim_{x \rightarrow \partial M} \tilde F(x) = l \in {\R} \cup \{ \infty \},
 \ x=(x_1,\dots , x_n),$$
if, for $\|x\| \rightarrow + \infty$ or $z \rightarrow z_0 \in
\partial M$, we have $\|\tilde F(x)\| \rightarrow l$,
where $\partial M$ denotes the boundary of $M\subset {\complex}^n$
and $\tilde F: \tilde M \rightarrow {\R}$,
$\tilde M$ given by (\ref{tildem}),  is the continuous map such that
$$F(z_1, \dots ,z_n)=\tilde F(x_1, \dots ,x_n), \ x_j=|z_j|^2 .$$}
\end{equation}
\end{theo}

Observe that  the maps $\Psi_0$, $\Psi_{hyp}$ and $\Psi_{FS}$ can
be described  explicitely (see (\ref{equa}), (\ref{equb}) and
(\ref{equc}) below). This is a rare phenomenon. In fact the proofs
of Theorems A, B, C, D and E above are existential and the
explicit form of the symplectic embedding or symplectomorphism
into a given complex space form is, in general, very hard to find.

\vskip 0.1cm

Theorem \ref{cor3} is an extension and a generalization of the
results  obtained by the first author and Fabrizio Cuccu in
\cite{alcu} for complete Reinhardt domains in ${\complex}^2$.
 Actually, all  the results obtained there
become a straightforward corollary of our Theorem \ref{cor3} (see
Example \ref{esempiohart} in Section \ref{examples}).

\vskip 0.1cm

Our second  result is Theorem \ref{teor2} below where we describe
geometric conditions on $\Phi$, related to Calabi's work on \K\
immersions, which implies  the existence of a special symplectic
immersion of $(M, \omega_{\Phi})$ in $({\real}^{2n}, \omega_0)$,
$n = \dim_{{\C}}M$ (and in particular the existence  of global
symplectic coordinates of $(M, \omega_{\Phi})$).

\begin{theo}\label{teor2}
Let $M \subseteq {\C}^n$ be a complex domain such that $0\in M$
endowed with a rotation invariant \K\ form $\omega_\Phi$. Assume
that  there exists a \K\  (i.e. a  holomorphic and isometric)
immersion  of $(M, g_\Phi)$ into some finite or infinite
dimensional complex space form, where $g_{\Phi}$ is the metric
associated to $\omega_{\Phi}$. Then, (\ref{genconda}) is satisfied
and hence there exists a special symplectic immersion $\Psi_0$ of
$(M, \omega_\Phi)$ into $({\C}^n, \omega_0)$, which is a global
symplectomorphism if and only if $\lim_{x \rightarrow
\partial M} \sum_{j=1}^n \frac{\partial{\tilde \Phi}}{\partial x_j} x_j =
+ \infty$. If $\sum_{j=1}^n \frac{\partial{\tilde \Phi}}{\partial
x_j} x_j < 1$ then there exists a symplectic immersion $\Psi_{FS}$
of $(M, \omega_\Phi)$ into $({\C}P^n, \omega_{FS})$ which is an
embedding if and only if $\lim_{x \rightarrow
\partial M} \sum_{j=1}^n \frac{\partial{\tilde \Phi}}{\partial x_j} x_j = 1$.
\end{theo}

\vskip 0.5cm The paper is organized as follows. In the next
section we prove Theorem \ref{cor3} and  Theorem \ref{teor2}. The
later  will follow  by an application of Calabi's results, which
will be briefly recalled in that section. Finally, in  Section
\ref{examples} we apply Theorem \ref{cor3}  to some important
cases. In particular we recover the results proved in \cite{alcu}
and we prove that each of the Ricci flat (but not flat) \K\ forms
on ${\complex}^2$ constructed by LeBrun  in \cite{lebrun1} admits
explicitely computable global symplectic coordinates. Observe that
this last result cannot be obtained by Theorem F above (see
Remark \ref{remarbefore} below).

\section{Proof of the main results}\label{sectionproofmain}
The following general lemma,  used in the proof of our main
results Theorem \ref{cor3} and  Theorem \ref{teor2}, describes the
structure of a special symplectic immersion between two complex
domains $M\subset{\complex}^n$ and $S\subset {\complex}^n$ endowed
with rotation invariant \K\ forms $\omega_{\Phi}$ and
$\omega_{\Xi}$ respectively. In all the paper we consider
smooth   maps
from $M$ into $S$ of the form
\begin{equation}\label{rotinv}
\Psi:M\rightarrow S, z \mapsto (\Psi_1(z)=\tilde \psi_1(x)
z_1,\dots , \Psi_n(z)=\tilde\psi_n(x)z_n),
\end{equation}
$z=(z_1,...,z_n)$, $x=(x_1,\dots ,x_n), x_j=|z_j|^2$
for some real functions $\tilde \psi_j: {\tilde M}
\rightarrow {\R}$, $j=1,\dots ,n$, where $\tilde M\subset {\real}^n$
is given by (\ref{tildem}).
A smooth map like (\ref{rotinv}) will be
called a  {\em special} map.

\begin{lemma}\label{lemmaa}
Let $M\subseteq {\C}^n$  and $S\subseteq {\C}^n$ be
complex domains
as above.
A special
map
$\Psi: M \rightarrow S , \; z
\mapsto (\Psi_1(z), \dots , \Psi_n(z)), $
is symplectic, namely   $\Psi^*(\omega_{\Xi}) = \omega_\Phi$,
 if and only if there exist constants $c_k \in \R$ such that
 the following equalities hold on $\tilde M$:
\begin{equation}\label{necsuff}
\tilde \psi^2_k \frac{\partial \tilde\Xi}{\partial x_k}(\Psi) =
\frac{\partial \tilde\Phi}{\partial x_k} + \frac{c_k}{x_k}, \; \;
\; k=1,\dots ,n,
\end{equation}
where $\tilde\Phi$ (resp. $\tilde\Xi$) is the function associated
to $\omega_{\Phi}$ (resp. $\omega_{\Xi}$),
and
$$ \frac{\partial \tilde\Xi}{\partial x_k}(\Psi)=
 \frac{\partial \tilde\Xi}{\partial x_k}
(\tilde\psi^2_1 x_1,\dots , \tilde \psi^2_n x_n), \; \; \;
k=1,\dots ,n.$$
\end{lemma}
\dimostr
From
$$\omega_{\Xi} = \frac{i}{2} \sum_{i,j = 1}^n \left(
\frac{\partial^2 \tilde\Xi}{\partial x_i
\partial x_j} \bar{z_j} z_i + \frac{\partial \tilde\Xi}{\partial
x_i} \delta_{ij} \right)_{x_1=|z_1|^2,\dots , x_n=|z_n|^2} d z_j \wedge d \bar z_i$$

one  gets

$$\Psi^*(\omega_{\Xi}) =
\frac{i}{2} \sum_{i,j=1}^n \left( \frac{\partial^2 \tilde
\Xi}{\partial x_i
\partial x_j}(\Psi) \Psi_i \bar\Psi_j + \frac{\partial \tilde\Xi}{\partial
x_j} (\Psi)\delta_{ij} \right)_{x_1=|z_1|^2,\dots , x_n=|z_n|^2}
 d\Psi_j \wedge d \bar\Psi_i ,$$
where
$$\frac{\partial^2 \tilde
\Xi}{\partial x_i
\partial x_j}(\Psi)=
\frac{\partial^2 \tilde
\Xi}{\partial x_i
\partial x_j}(\tilde\psi^2_1 x_1,\dots , \tilde \psi^2_n x_n).$$

If one denotes by
$$\Psi^*(\omega_{\Xi})=
\Psi^*(\omega_{\Xi})_{(2, 0)}+
\Psi^*(\omega_{\Xi})_{(1, 1)}+
\Psi^*(\omega_{\Xi})_{(0, 2)}$$
the decomposition of
$\Psi^*(\omega_{\Xi})$
into  addenda of type
$(2,0), (1,1)$ and $(0,2)$
one has:

\begin{equation}\label{02}
\Psi^*(\omega_{\Xi})_{(2, 0)}=
 \frac{i}{2} \sum_{i,j,k,l=1}^n \left( \frac{\partial^2 \tilde
\Xi}{\partial x_i
\partial x_j}(\Psi) \Psi_i \bar{\Psi_j} + \frac{\partial \tilde\Xi}{\partial
x_j} (\Psi)\delta_{ij} \right) \frac{\partial \Psi_j}{\partial z_k}
\frac{\partial \bar{\Psi_i}}{\partial z_l}  dz_k \wedge dz_l
\end{equation}

\begin{equation}\label{11}
\Psi^*(\omega_{\Xi})_{(1, 1)}
= \frac{i}{2} \sum_{i,j,k,l=1}^n\left( \frac{\partial^2 \tilde
\Xi}{\partial x_i
\partial x_j}(\Psi) \Psi_i \bar{\Psi_j} + \frac{\partial \tilde\Xi}{\partial
x_j} (\Psi)\delta_{ij} \right)  \left( \frac{\partial \Psi_j}{\partial
z_k} \frac{\partial \bar{\Psi_i}}{\partial \bar{z_l}}  -
\frac{\partial \Psi_j}{\partial \bar{z_l}} \frac{\partial
\bar{\Psi_i}}{\partial z_k} \right) dz_k \wedge d \bar{z_l}
\end{equation}

\begin{equation}\label{20}
\Psi^*(\omega_{\Xi})_{(0, 2)}=
\frac{i}{2} \sum_{i,j,k,l=1}^n\left( \frac{\partial^2 \tilde
\Xi}{\partial x_i
\partial x_j}(\Psi) \Psi_i \bar{\Psi_j} + \frac{\partial \tilde\Xi}{\partial
x_j} (\Psi)\delta_{ij} \right)  \; \frac{\partial \Psi_j}{\partial
\bar{z_k}} \frac{\partial \bar{\Psi_i}}{\partial \bar{z_l}} \; d
\bar{z_k} \wedge d \bar{z_l}.
\end{equation}
(Here and  below,   with a slight abuse of notation,  we are
omitting the fact that all the previous expressions have to be
evaluated at $x_1=|z_1|^2,\dots,  x_n=|z_n|^2$. ) Since $\Psi_j
(z)=  \tilde \psi_j( |z_1|^2,...,|z_n|^2)z_j$,
 one has:

\begin{equation}\label{ins1}
\frac{\partial \Psi_i}{\partial z_k} = \frac{\partial
\tilde \psi_i}{\partial x_k} z_i \bar{z_k} + \tilde \psi_i
\delta_{ik},\
\frac{\partial \Psi_i}{\partial \bar{z_k}} =
\frac{\partial \tilde \psi_i}{\partial x_k} z_k z_i
\end{equation}
and
\begin{equation}\label{ins2}
\frac{\partial \bar\Psi_i}{\partial \bar{z_k}} =
\frac{\partial \tilde \psi_i}{\partial x_k} z_k \bar{z_i} + \tilde
\psi_i \delta_{ik},\
 \frac{\partial \bar\Psi_i}{\partial z_k} =
\frac{\partial \tilde \psi_i}{\partial x_k} \bar{z_k} \bar{z_i},
\end{equation}
By inserting (\ref{ins1}) and (\ref{ins2}) into (\ref{02})
and (\ref{11})
after a long, but straightforward
computation,
 one obtains:
\begin{equation}\label{20A}
\Psi^*(\omega_{\Xi})_{(2, 0)}=
 \frac{i}{2} \sum_{k,l=1}^n
\frac{A_{kl}}{2} \bar z_k \bar z_l \; dz_k\wedge dz_l
\end{equation}

and

\begin{equation}\label{11A}
\Psi^*(\omega_{\Xi})_{(1, 1)}
= \frac{i}{2}
 \sum_{k,l=1}^n
\left[ (\frac{A_{kl}+A_{lk}}{2}+
\frac{\partial^2 \tilde
\Xi}{\partial x_k
\partial x_l} (\Psi)\tilde\psi_k^2 \tilde\psi_l^2)\bar z_kz_l+
\frac{\partial \tilde \Xi}{\partial x_k} (\Psi)\delta_{kl}
\tilde\psi_k^2\right] dz_k\wedge d\bar z_l,
\end{equation}

where

\begin{equation}\label{Akl}
A_{kl}=\frac{\partial \tilde\Xi}{\partial x_k} (\Psi)\frac{\partial \tilde
\psi_k^2}{\partial x_l} + \tilde \psi_k^2 \sum_{j=1}^n
\frac{\partial^2 \tilde\Xi}{\partial x_j
\partial x_k} (\Psi)\frac{\partial \tilde \psi_j^2}{\partial x_l} |z_j|^2.
 \end{equation}

Now, we assume that
$$\Psi^*(\omega_{\Xi}) = \omega_\Phi
= \frac{i}{2} \sum_{k,l = 1}^n \left( \frac{\partial^2 \tilde
\Phi}{\partial x_k
\partial x_l} \bar{z_k} z_l + \frac{\partial \tilde\Phi}{\partial
x_l} \delta_{lk} \right)_{x_1=|z_1|^2,\dots , x_n=|z_n|^2}  d z_k \wedge d \bar{z_l}.$$

Then the terms  $\Psi^*(\omega_{\Xi})_{(2, 0)}$
and $\Psi^*(\omega_{\Xi})_{(0, 2)}$ are
equal to  zero.
This is equivalent to the fact
that (\ref{Akl})
is symmetric in $k,l$.

Hence, by setting
$\Gamma_l =  \tilde \psi^2_l \; \frac{\partial
\tilde\Xi}{\partial x_l}(\Psi),\  l=1,\dots ,n$
equation
(\ref{11A}) becomes
$$\Psi^*(\omega_{\Xi})_{(1, 1)}
= \frac{i}{2}
 \sum_{k,l=1}^n
\left[ (A_{kl}+
\frac{\partial^2 \tilde
\Xi}{\partial x_k
\partial x_l}(\Psi) \tilde\psi_k^2 \tilde\psi_l^2)\bar z_kz_l+
\frac{\partial \tilde \Xi}{\partial x_k}(\Psi)\delta_{kl}
\tilde\psi_k^2\right] dz_k\wedge d\bar z_l =
$$
\begin{equation}\label{symp}
=\frac{i}{2} \sum_{k,l = 1}^n \left( \frac{\partial
\Gamma_l}{\partial x_k} \bar{z_k} z_l + \Gamma_k \delta_{kl}
\right) d z_k \wedge d \bar{z_l} .
\end{equation}

So, $\Psi^*(\omega_{\Xi}) = \omega_\Phi$ implies

$$\frac{i}{2} \sum_{k,l = 1}^n \left( \frac{\partial
\Gamma_l}{\partial x_k} \bar{z_k} z_l + \Gamma_k \delta_{lk}
\right) d z_k \wedge d \bar{z_l}  = \frac{i}{2} \sum_{k,l = 1}^n
\left( \frac{\partial^2 \tilde\Phi}{\partial x_k
\partial x_l} \bar{z_k} z_l + \frac{\partial \tilde\Phi}{\partial
x_l} \delta_{kl} \right) d z_k \wedge d \bar{z_l}.$$

In this equality, we distinguish the cases $l \neq k$ and $l = k$
and  get respectively

$$\frac{\partial \Gamma_l}{\partial x_k} = \frac{\partial^2 \tilde
\Phi}{\partial x_k \partial x_l} \; \; \;(k \neq l)$$

and

$$\frac{\partial \Gamma_k}{\partial x_k} x_k + \Gamma_k =
\frac{\partial^2 \tilde\Phi}{\partial x_k^2} x_k + \frac{\partial
\tilde\Phi}{\partial x_k} .$$

By defining $A_k = \Gamma_k - \frac{\partial \tilde
\Phi}{\partial x_k}$, these equations become respectively

$$\frac{\partial A_k}{\partial x_l} = 0 \; \; \; (l\neq k)$$

and

$$\frac{\partial A_k}{\partial x_k} \  x_k = - A_k .$$

The first  equation implies that $A_k$
does not depend on $x_l$
and so by the second one
we have
\begin{equation}\label{Ak}
A_k = \Gamma_k - \frac{\partial \tilde \Phi}{\partial x_k} =
\frac{c_k}{x_k},
\end{equation}
for some constant $c_k \in {\real}$, i.e.
$$\Gamma_k = \tilde \psi^2_k \frac{\partial
\tilde\Xi}{\partial x_k}(\Psi) = \frac{\partial
\tilde\Phi}{\partial x_k} + \frac{c_k}{x_k}, \; \; \; \ k=1,\dots
,n,$$ namely (\ref{necsuff}).

 In order to prove the converse of Lemma \ref{lemmaa},
notice that by
differentiating (\ref{necsuff}) with respect to $l$ one gets:

$$\frac{\partial^2 \tilde\Phi}{\partial x_k \partial x_l} - \frac{c_k}{x_k^2} \delta_{kl} =
A_{kl}+
 \frac{\partial^2
\tilde\Xi}{\partial x_k
\partial x_l} \tilde \psi_k^2 \tilde \psi_l^2$$
with $A_{kl}$ given by (\ref{Akl}).
By $\frac{\partial^2 \tilde\Phi}{\partial x_k
\partial x_l} = \frac{\partial^2 \tilde\Phi}{\partial x_l \partial
x_k}$ and
$ \frac{\partial^2
\tilde\Xi}{\partial x_k
\partial x_l} \tilde \psi_k^2 \tilde \psi_l^2=
 \frac{\partial^2
\tilde\Xi}{\partial x_l
\partial x_k} \tilde \psi_l^2 \tilde \psi_k^2$  one gets
$A_{kl}=A_{lk}$. Then, by (\ref{20A}),  the addenda of type (2,0)
(and (0,2)) in $\Psi^*(\omega_{\Xi})$ vanish. Moreover, by
(\ref{Akl}) and (\ref{symp}),
 it follows that  $\Psi^*(\omega_{\Xi})=\omega_\Phi$.
\fdim

\vskip 0.3cm

 In the proof of  Theorem \ref{cor3}
 we also  need the following lemma whose proof
follows by Lemma
\ref{lemmaa}, or by a direct computation.
\begin{lemma}\label{lemmadual}
The map
$f: {\C}H^n \rightarrow
{\C}^n$
given by
\begin{equation}\label{eqfh}
(z_1, \dots, z_n)\mapsto \left( \frac{z_1}{\sqrt{1 - \sum_{i=1}^n |z_i|^2 }}, \dots ,\frac{z_n}{\sqrt{1 - \sum_{i=1}^n |z_i|^2 }}
\right)
\end{equation}
is a special global diffeomorphism
satisfying
\begin{equation}\label{s1}
f^*(\omega_0)=\omega_{hyp}
\end{equation}
and
\begin{equation}\label{s2}
f^*(\omega_{FS})=\omega_{0},
\end{equation}
where, in the second equation,
we are looking at ${\C}^n \stackrel{i}{\hookrightarrow}
{\complex}P^n$ as the affine chart $Z_0\neq 0$ in ${\complex}P^n$
endowed with the restriction of the Fubini--Study form $\omega_{FS}$
and where   $\omega_0$
denotes the restriction of the flat form  of ${\complex}^n$
to ${\complex}H^n\subset {\complex}^n$.
\end{lemma}

We are now in the position to prove our first result.

\vskip 0.5cm

\noindent
{\bf Proof of Theorem \ref{cor3}}
First of all observe that under
assumption (\ref{intersezassi}),
the $c_k$'s appearing in the
statement of Lemma \ref{lemmaa} are forced to be zero.
 So, the existence of
a special symplectic immersion $\Psi: M \rightarrow S$ is
equivalent to
\begin{equation}\label{remarcon}
\tilde \psi^2_k \frac{\partial \tilde\Xi}{\partial x_k}(\Psi) =
\frac{\partial \tilde\Phi}{\partial x_k}, \  k=1, \dots, n .
\end{equation}
If we further assume  $(S={\complex}^n,  \omega_{\Xi} =\omega_0)$,
namely $\tilde\Xi =\sum_{j=1}^nx_j$, then condition
 (\ref{necsuff}) reduces to
$$\tilde \psi^2_k= \frac{\partial
\tilde\Phi}{\partial x_k}, \  k=1, \dots ,n ,$$ and hence
(\ref{cond0})  follows by Lemma \ref{lemmaa}. Further  $\Psi_0$ is
given  by:
 \begin{equation}\label{equa}
\Psi_0(z) =
\left(\sqrt{\frac{\partial \tilde \Phi}{\partial x_1}}\  z_1,\cdots, \sqrt{\frac{\partial \tilde \Phi}{\partial x_n}} \  z_n
\right) _{x_i = |z_i|^2}
\end{equation}
In order to prove (i) when
$(S={\complex}H^n,  \omega_{\Xi} =\omega_{hyp})$
observe  that since the composition of two special maps
is a special map it follows  by  (\ref{s1}) that
the existence of a  special symplectic map
$\Psi: (M, \omega_{\Phi})\rightarrow({\complex}H^n,  \omega_{hyp})$
gives rise to a  special symplectic map
$f\circ \Psi: (M, \omega_{\Phi})\rightarrow({\complex}^n,  \omega_{0})$.
The later is uniquely determined by  the previous case, i.e.
$\Psi_0=f\circ \Psi$.  So $\Psi_{hyp}:=\Psi =f^{-1}\circ\Psi_0$
and  since   the inverse of $f$ is given by
$$f^{-1}: {\C}^n \rightarrow
{\C}H^n,
z\mapsto \left( \frac{z_1}{\sqrt{1 + \sum_{i=1}^n |z_i|^2 }}, \dots ,\frac{z_n}{\sqrt{1 + \sum_{i=1}^n |z_i|^2 }}
\right),$$
one obtains:
\begin{equation}\label{equb}
\Psi_{hyp}(z) = \left( \sqrt{ \frac{ \frac{\partial
\tilde\Phi}{\partial x_1} }{ 1+ \sum_{k=1}^n \frac{\partial
\tilde\Phi}{\partial x_k} x_k} \ } z_1,\dots , \sqrt{ \frac{
\frac{\partial \tilde\Phi}{\partial x_n} }{ 1+ \sum_{k=1}^n
\frac{\partial \tilde\Phi}{\partial x_k} x_k} }\  z_n \right)
_{x_i = |z_i|^2}.
\end{equation}

In order to prove (ii),
 notice that by (\ref{s2})
a  special symplectic map
$\Psi: (M, \omega_{\Phi})\rightarrow({\complex}^n,  \omega_{FS})$
is uniquely determined by  the  special symplectic map
$\Psi_0=f^{-1}\circ \Psi: (M, \omega_{\Phi})\rightarrow({\complex}H^n,  \omega_{0})\subset
({\complex}^n,  \omega_{0})$
and therefore (\ref{conda})  is  a straightforward
consequence of the previous case (i).
 Furthermore $\Psi_{FS}$
 is given by
\begin{equation}\label{equc}
\Psi_{FS}(z) = \left( \sqrt{ \frac{ \frac{\partial
\tilde\Phi}{\partial x_1} }{ 1 - \sum_{k=1}^n \frac{\partial
\tilde\Phi}{\partial x_k} x_k} }\  z_1,\dots ,\sqrt{ \frac{
\frac{\partial \tilde\Phi}{\partial x_n} }{ 1 - \sum_{k=1}^n
\frac{\partial \tilde\Phi}{\partial x_k} x_k} }\  z_n \right)
_{x_i = |z_i|^2}
\end{equation}

Finally, notice that conditions  (\ref{genconda}) and
(\ref{gencondb})
 for  the special map (\ref{equa})
(resp. (\ref{equc})) are  equivalent to $\Psi_0^{-1}(\{ 0 \}) = \{
0 \}$ (resp. $\Psi_{FS}^{-1}(\{ 0 \}) = \{ 0 \}$) and to the
properness of $\Psi_0$ (resp. $\Psi_{FS}$). Hence,   the fact that
under these conditions
 $\Psi_0$ (resp. $\Psi_{FS}$) is a global diffeomorphism follows
 by standard topological arguments.
 \fdim

\begin{remar} \label{remarkarbsmall} \rm
Observe that, by Theorem \ref{cor3}, if $(M, \omega_{\Phi})$
admits a special symplectic immersion into $({\C}^n,
\omega_{FS})$, then it admits a special symplectic immersion in
$({\C}^n, \omega_{0})$ (or $({\C}H^n, \omega_{hyp})$). The
converse is false even if one restricts to an arbitrary small open
set $U \subseteq M$ endowed with the restriction of
$\omega_{\Phi}$ (see Remark \ref{remarex} below).
\end{remar}

In order to prove our second result  (Theorem \ref{teor2})
we  briefly recall
Calabi's work on \K\ immersions and his fundamental  Theorem \ref{criterium}.
 We refer the reader to  \cite{ca}
 for details and further results (see also  \cite{hartogs} and \cite{diastherm}).

\vskip 0.5cm

\noindent
{\bf Calabi's work}
In his seminal paper Calabi \cite{ca} gave a complete answer to the problem of
the existence and uniqueness of \K\
immersions of a \K\ manifold $(M, g)$ into a finite or
infinite dimensional complex space form.
Calabi's first observation was that if  such  \K\ immersion exists
then the metric $g$ is forced to be real analytic being the
pull-back via a holomorphic map of the real analytic metric  of a
complex space form. Then in a neighborhood of every point $p\in
M$, one can introduce a very special \K\ potential $D^g_p$ for the
metric $g$, which Calabi christened {\em diastasis}. The
construction goes as follows.
Take a a real-analytic \K\ potential  $\Phi$ around the point $p$
(it exists since $g$ is real analytic). By duplicating the
variables $z$ and $\bar z$ $\Phi$ can be complex analytically
continued to a function $\hat\Phi$ defined in a neighborhood $U$
of the diagonal containing $(p, \bar p)\in M\times\bar M$ (here
$\bar M$ denotes the manifold conjugated to $M$). The {\em
diastasis function} is the \K\ potential $D^g_p$ around $p$
defined by
$$D^g_p(q)=\hat\Phi (q, \bar q)+
\hat\Phi (p, \bar p)-\hat\Phi (p, \bar q)- \hat\Phi (q, \bar p).$$

\begin{Example}\label{flat}\normalshape
Let $g_0$ be the Euclidean metric on ${\complex}^N$,
$N\leq\infty$, namely the metric whose associated \K\ form is
given by $\omega_{0}=\frac{i}{2}\sum_{j=1}^Ndz_j\wedge d\bar z_j$.
Here ${\complex}^{\infty}$ is the complex  Hilbert space
$l^2(\complex)$ consisting of sequences $(z_j)_{j \geq 1}$, $z_j
\in {\complex}$ such that $\sum_{j=1}^{+\infty}|z_j|^2<+\infty$.
The diastasis function $D_0^{g_{0}}:{\complex}^N\rightarrow
{\real}$ around the origin $0\in {\complex}^N$ is given by
\begin{equation}\label{equflat}
D_0^{g_{0}}(z)=\sum_{j=1}^N|z_j|^2.
\end{equation}
\end{Example}
\begin{Example}\label{proj}\normalshape
Let $(Z_0, Z_1, \dots ,Z_N)$ be the homogeneous coordinates in the
complex projective space in ${\complex}P^N$, $N\leq\infty$,
endowed with the Fubini--Study metric $g_{FS}$. Let $p=[1, 0,\dots
,0]$. In the affine chart $U_0=\{Z_0\neq 0\}$ endowed with
coordinates $(z_1,\dots ,z_N), z_j=\frac{Z_j}{Z_0}$ the diastasis
around $p$ reads as:
\begin{equation}\label{equproj}
D^{g_{FS}}_{p}(z)= \log (1+\sum _{j=1}^{N}|z_j|^2).
\end{equation}
\end{Example}
\begin{Example}\label{hyp}\normalshape
Let ${\complex}H^N=\{z\in {\C}^N|\ \sum
_{j=1}^{N}|z_j|^2<1\}\subset {\complex}^N, N\leq\infty$ be the
complex hyperbolic space endowed with the hyperbolic metric
$g_{hyp}$. Then the diastasis around the origin is given by:
\begin{equation}\label{equhyp}
D^{g_{hyp}}_{0}(z)= -\log (1-\sum _{j=1}^{N}|z_j|^2).
\end{equation}
\end{Example}

A very useful characterization of the diastasis (see below) can be
obtained as follows. Let $(z)$ be a system of complex coordinates
in a neighbourhood of $p$ where $D_p^g$ is defined. Consider its
power series development:
\begin{equation}\label{diastflat}
D_p^g(z)=\sum _{j, k\geq 0} a_{jk}(g)z^{m_j}\bar z^{m_k},
\end{equation}

\noindent where we are using the following convention: we arrange
every $n$-tuple of nonnegative integers as the sequence
$$m_{j}=(m_{1, j}, m_{2, j},\dots , m_{n, j}) _{j=0, 1,\dots}$$
such that $m_0=(0, \dots , 0)$, $|m_{j}|\leq |m_{j+1}|$, with
$|m_{j}|=\sum_{\alpha =1}^n m_{\alpha , j}$ and
$z^{m_{j}}=\prod_{\alpha =1}^{n} (z_{\alpha})^{m_{\alpha , j}}$.
Further, we order all the $m_j$'s with the same $|m_j|$ using the
lexicographic order.

Characterization of the diastasis: {\em Among all the
potentials the diastasis is characterized by the fact that in
every coordinates system $(z)$ centered in $p$ the coefficients
$a_{jk}(g)$ of the expansion (\ref{diastflat}) satisfy $a_{j
0}(g)=a_{0 j}(g)=0$ for every nonnegative integer $j$.}

\begin{defin}\label{full}
A \K\ immersion $\varphi$ of $(M, g)$ into a complex space form
$(S, G)$ is said to be full if  $\varphi (M)$ is not contained in
a proper complex totally geodesic submanifold of $(S, G)$.
\end{defin}

\begin{defin}\label{res}
Let  $g$ be a real analytic  \K\ metric on a complex manifold $M$.
The metric $g$ is said to be resolvable of rank  $N$ if the
$\infty\times\infty$ matrix $a_{jk}(g)$ given by
({\ref{diastflat}}) is positive semidefinite and of rank $N$.
Consider the function $e^{D_p^g}-1$ (resp. $1-e^{-D_p^g}$ ) and
its power series development:
\begin{equation}\label{diastproj}
e^{D_p^g}-1= \sum _{j, k\geq 0}b_{jk}(g)z^{m_j} \bar{z}^{m_k}.
\end{equation}
(resp.
\begin{equation}\label{diasthyp}
1-e^{-D_p^g}= \sum _{j, k\geq 0}c_{jk}(g)z^{m_j} \bar{z}^{m_k}.)
\end{equation}

The metric $g$ is said to be $1$-resolvable (resp.
$-1$-resolvable) of rank $N$ at $p$ if the $\infty\times\infty$
matrix $b_{jk}(g)$ (resp.  $c_{jk}(g)$) is positive semidefinite
and of  rank $N$.
\end{defin}

We are now in the position to state Calabi's fundamental theorem
and to  prove our  Theorem \ref{teor2}.

\begin{theo}{\bf (Calabi)}\label{criterium}
Let $M$ be a complex manifold endowed with a real analytic \K\ \
metric $g$. A neighbourhood of a point $p$ admits a (full) \K\ \
immersion into  $({\C}^N, g_0)$ if and only if $g$ is resolvable
of rank at most (exactly) $N$ at  $p$. A neighbourhood of a point
$p$ admits a (full) \K\ immersion into  $({\C}P^N, g_{FS})$ (resp.
$({\C}H^N, g_{hyp})$ ) if and only if $g$ is $1$-resolvable (resp.
$-1$-resolvable) of rank at most (exactly) $N$ at  $p$.
\end{theo}

\vskip 0.3cm

\noindent {\bf Proof of Theorem \ref{teor2} } Without loss of
generality we can assume $\Phi (0)=0$. Then, it follows by the
characterization
 of the diastasis function
that $\Phi$ is indeed the (globally defined) diastasis function
for the \K\ metric $g_\Phi$ (associated to $\omega_\Phi$) around
the origin,
 namely $\Phi=D_0^{g_\Phi}$.
 Since $\Phi=D_0^{g_\Phi}$ is rotation invariant, namely it
depends only on $|z_1|^2,\dots , |z_n|^2$, the matrices
$a_{jk}(g)$,  $b_{jk}(g)$ and $c_{jk}(g)$ above   are diagonal,
i.e.
\begin{equation}\label{bj}
a_{jk}(g)=a_j\delta_{jk}, \ b_{jk}(g)=b_j\delta_{jk}, \
c_{jk}(g)=c_j\delta_{jk}, \ a_j, b_j, c_j \in {\real}.
\end{equation}
Therefore, by Calabi's Theorem \ref{criterium} if $(M, g_\Phi)$
 admits a (full) \K\ immersion into
 $({\C}^N, g_0)$
(resp.  $({\C}P^N, g_{FS})$ or
 $({\C}H^N, g_{hyp})$) then
 all the $a_j's$ (resp. the $b_j's$ or
the  $c_j's$) are greater or equal than $0$ and at most (exactly)
$N$ of them are positive.
Moreover,
 it follows by the fact
 that the metric $g_\Phi$
 is positive definite (at $0\in M$)
 that  the  coefficients
$a_k$ (resp. $b_k$ or $c_k$), $k= 1, \dots, n$,  are strictly
greater than zero . Hence,
 by using (\ref{diastflat})
(resp. (\ref{diastproj}) or (\ref{diasthyp})) with $p=0$ and
$g=g_{\Phi}$ we get $ \frac{\partial \tilde \Phi}{\partial x_k}(x)
= a_k + P_0(x)$ (resp. $\frac{\partial \tilde \Phi}{\partial
x_k}(x) = \frac{b_k + P_+(x)}{1 + \sum_j b_j x^{m_j}}$ or
$\frac{\partial \tilde \Phi}{\partial x_k}(x) = \frac{c_k +
P_-(x)}{1 - \sum_j c_j x^{m_j}})$ where $P_0$ (resp. $P_+$ or
$P_-$) is a polynomial with non-negative coefficients
in the variables $x=(x_1, \dots ,x_n), x_j=|z_j|^2$.
Hence
condition (\ref{cond0}) above is satisfied. The last two
assertions of Theorem \ref{teor2} are immediately
 consequences of Theorem \ref{cor3}. \fdim

\section{Applications and further results}\label{examples}

\begin{Example}\label{esempiohart}\rm
(cfr. \cite{alcu}) Let $x_0 \in \R^+ \cup \{ + \infty \}$ and let
$F: [0, x_0) \rightarrow (0, + \infty)$ be a non-increasing smooth
function. Consider the domain

$$D_F = \{ (z_1, z_2) \in {\C}^2 \; | \; |z_1|^2 < x_0, \ |z_2|^2 < F(|z_1|^2)
\}$$

\noindent endowed with the 2-form $\omega_F = \frac{i}{2} \partial
\overline{\partial} \log \frac{1}{F(|z_1|^2) - |z_2|^2}$. If the
function $A(x) = - \frac{x F'(x)}{F(x)}$ satisfies $A'(x) >0$ for
every $x \in [0, x_0)$, then $\omega_F$ is a \K\ form on $D_F$ and
$(D_F, \omega_F)$ is called the {\em complete Reinhardt domain}
 associated with $F$. Notice that
$\omega_F$ is rotation invariant with associated real function
$\tilde F(x_1, x_2) = \log \frac{1}{F(x_1) - x_2}$. We now apply
Theorem \ref{cor3} to $(D_F, \omega_F)$. We have

$$\frac{\partial \tilde F}{\partial x_1} = - \frac{F'(x_1)}{F(x_1) - x_2}>0, \; \; \; \frac{\partial \tilde F}{\partial x_2} = \frac{1}{F(x_1) - x_2}>0,\  x_j=|z_j|^2, \  j=1, 2.$$

So, by Theorem \ref{cor3}, $(D_F, \omega_F)$ admits a
special symplectic immersion in $({\C}^2, \omega_0)$ (and in
$({\C}H^2, \omega_{hyp})$). Moreover, this immersion is a global
symplectomorphism only when
$$\frac{\partial \tilde
F}{\partial x_1} x_1 +
\frac{\partial \tilde
F}{\partial x_2} x_2
 = \frac{x_2 - F'(x_1)x_1}{F(x_1) - x_2}$$
tends to infinity on the boundary of $D_F$. For example, let $F:
[0, +\infty) \rightarrow {\R}^+$ given by $F(x) = \frac{c}{c+x}$,
$c>0$ (resp.  $F(x) = \frac{1}{(1+x)^p}$, $p \in {\N}^+$). Then
$\sum_{i=1}^2 \frac{\partial \tilde F}{\partial x_i} x_i =
\frac{x_2(c+x_1)^2 + c x_1}{c(c+x_1) - x_2(c + x_1)^2}$ (resp.
$\sum_{i=1}^2 \frac{\partial \tilde F}{\partial x_i} x_i =
\frac{x_2 + p x_1(1+x_1)^{-p-1}}{(1+x_1)^{-p} - x_2}$) does not
tend to infinity, for $t \rightarrow \infty$, along the curve
$x_1=t$, $x_2 = \frac{\varepsilon c}{c+t}$, for any $\varepsilon
\in (0,1)$ (resp. does not tend to infinity, for $t \rightarrow
\infty$, along the curve $x_1=t$, $x_2 = \varepsilon (1+t)^{-p}$,
for any $\varepsilon \in (0,1)$). On the other hand, one verifies
by a straight calculation that, if $F: [0, +\infty) \rightarrow
{\R}^+$ is given by $F(x) = e^{-x}$ (resp. $F: [0, 1) \rightarrow
{\R}^+$, $F(x) = (1-x)^p$, $p>0$), then $\sum_{i=1}^2
\frac{\partial \tilde F}{\partial x_i} x_i = \frac{x_2 + e^{-x_1}
x_1}{e^{-x_1} - x_2}$ (resp. $\sum_{i=1}^2 \frac{\partial \tilde
F}{\partial x_i} x_i = \frac{x_2 + p x_1 (1- x_1)^{p-1}}{(1-x_1)^p
- x_2}$) tends to infinity on the boundary of $D_F$. We then
recover the conclusions of Examples  3.3, 3.4, 3.5, 3.6 in
\cite{alcu}.

\end{Example}

\begin{Example}\rm\label{esempioloc1}
Let us endow ${\C}^2 \setminus \{ 0 \}$ with the rotation
invariant \K\ form $\omega_{\Phi} = \frac{i}{2} \partial \bar
\partial \Phi$
with associated real function
$$\tilde\Phi (x_1, x_2) = a \log(x_1 + x_2) + b(x_1 + x_2) + c , \  a,b,c >0 .$$
The metric $g_{\Phi}$ associated to $\omega_{\Phi}$
is used in  \cite{simanca} (see also   \cite{lebrun2})
for the construction of
\K\ metrics of constant scalar curvature on bundles on ${\complex}P^{n-1}$.

Since $\frac{\partial \tilde \Phi}{\partial x_i} = b +
\frac{a}{x_1+x_2}
>0$, by Theorem \ref{cor3} there exists a special symplectic
immersion of $({\C}^2 \setminus \{ 0 \}, \omega_\Phi)$ in
$({\C}^2, \omega_0)$ (or in $({\C}H^2, \omega_{hyp})$).
\end{Example}

\begin{Example}\label{esempioloc2}\rm
Let us endow ${\C}^2 \setminus \{ 0 \}$ with the metric
$\omega_{\Phi} = \frac{i}{2} \partial \bar \partial \Phi$, where
$$\tilde \Phi = \sqrt{r^4 + 1} + 2 \log r - \log(\sqrt{r^4 + 1} +
1), \ r = \sqrt{|z_1|^2 + |z_2|^2} .$$

The metric $g_{\Phi}$ is used in \cite{joyce} for the construction
of the Eguchi--Hanson metric. A straight calculation shows that
$$\frac{\partial \tilde\Phi}{\partial x_i} = \frac{\partial
\tilde\Phi}{\partial r} \frac{\partial r}{\partial x_i} = \left[
\frac{4 r^3}{2 \sqrt{r^4+1}} \left( 1 - \frac{1}{\sqrt{r^4+1} +1}
\right) + \frac{2}{r} \right] \frac{1}{2 r}
> 0 ,$$ so by Theorem \ref{cor3} there exists a special symplectic immersion of \linebreak $({\C}^2 \setminus \{ 0 \},
\omega_\Phi)$ in $({\C}^2, \omega_0)$ (or in $({\C}H^2,
\omega_{hyp})$).

\end{Example}

\begin{remar}\label{remarex} \rm

Notice that in the previous Example \ref{esempioloc2} one has

$$\frac{\partial \tilde\Phi}{\partial x_1} x_1
+ \frac{\partial \tilde\Phi}{\partial x_2} x_2 =
\frac{r^4}{\sqrt{r^4+1}} \left( 1 - \frac{1}{\sqrt{r^4+1} +1}
\right) + 1 > 1 ,$$ so again by Theorem \ref{cor3} it does not
exist a special symplectic immersion of $({\C}^2 \setminus \{ 0
\}, \omega_\Phi)$ in $({\C}^2, \omega_{FS})$.
Moreover, such an immersion does not exist for any
arbitrarily small $U \subseteq {\C}^2 \setminus \{ 0 \}$ endowed
with the restriction of $\omega_{\Phi}$ (cfr.  Remark
\ref{remarkarbsmall} above).

\end{remar}

\begin{Example}\rm
 In \cite{lebrun1} Claude LeBrun constructed
the following  family of  \K\ forms on ${\C}^2$  defined by
$\omega_m = \frac{i}{2} \partial \bar \partial \Phi_m$, where
$$\Phi_m(u,v) = u^2 + v^2 + m (u^4 + v^4),\  m \geq 0$$
 and $u$ and $v$ are implicitly defined by
 $$|z_1| =
e^{m(u^2-v^2)} u,\  |z_2| = e^{m(v^2-u^2)} v .$$
 For $m = 0$ one
gets the flat metric, while for $m>0$ each of the metrics of this
family represents
 the first example of complete Ricci flat (non-flat) metric on ${\C}^2$ having the same
volume form of the flat metric $\omega_0$, namely
$\omega_m\wedge\omega_m=\omega_0\wedge\omega_0$.
 Moreover, for $m>0$,
these metrics are isometric (up to dilation and rescaling)  to the
Taub-NUT  metric.

Now, with the aid of Theorem \ref{cor3}, we  prove that for every
$m$ the \K\ manifold  $({\C}^2, \omega_{m})$ admits global symplectic
coordinates.
Set  $u^2 = U$, $v^2 = V$. Then $\tilde \Phi_m$ (the function
associated to $\Phi_m$) satisfies:

$$\frac{\partial\tilde\Phi_m}{\partial x_1} = \frac{\partial \tilde\Phi_m}{\partial
U}\frac{\partial U}{\partial x_1} + \frac{\partial
\tilde\Phi_m}{\partial V} \frac{\partial V}{\partial x_1} ,$$

$$\frac{\partial\tilde \Phi_m}{\partial x_2} = \frac{\partial\tilde \Phi_m}{\partial
U}\frac{\partial U}{\partial x_2} + \frac{\partial
\tilde\Phi_m}{\partial V} \frac{\partial V}{\partial x_2} ,$$
where $x_j=|z_j|^2, j=1, 2$.
In order to calculate $\frac{\partial U}{\partial x_j}$ and
$\frac{\partial V}{\partial x_j}, j=1, 2$,  let us consider the map
$$G:
{\R}^2 \rightarrow {\R}^2, \ (U, V) \mapsto  (x_1=e^{2m(U-V)} U,\
x_2=e^{2m(V-U)} V)$$ and its Jacobian matrix

\begin{displaymath}
J_G = \left( \begin{array}{cc}
(1 + 2mU) \; e^{2m(U-V)} &  - 2 m U \; e^{2m(U-V)} \\
- 2 m V \; e^{2m(V-U)} &  (1 + 2m V) \; e^{2m(V-U)}
\end{array} \right) .
\end{displaymath}

We have $det J_G = 1 + 2 m (U + V) \neq 0$, so

$$J_G^{-1} = J_{G^{-1}} = \frac{1}{1 + 2 m (U + V)} \left( \begin{array}{cc}
(1 + 2m V) e^{2m(V-U)} &  2 m U e^{2m(U-V)} \\
2 m V e^{2m(V-U)} &  (1 + 2m U) e^{2m(U-V)}
\end{array} \right) .$$

Since
$J_{G^{-1}} = \left( \begin{array}{cc}
\frac{\partial U}{\partial x_1} &  \frac{\partial U}{\partial x_2} \\
\frac{\partial V}{\partial x_1} &  \frac{\partial V}{\partial x_2}
\end{array} \right)$,
by a straightforward calculation we get
$$\frac{\partial\tilde\Phi_m}{\partial x_1} = (1 + 2m V) e^{2m(V-U)}  > 0, \ \frac{\partial\tilde\Phi_m}
{\partial x_2} = (1 + 2m U) e^{2m(U-V)}
> 0, $$
and
$$\lim_{\|x\| \rightarrow+ \infty} (\frac{\partial\tilde\Phi_m}{\partial x_1} x_1 + \frac{\partial\tilde\Phi_m}{\partial
x_2}x_2) = \lim_{\|x\|  \rightarrow +\infty}  (U + V + 4 m UV)=
+\infty, $$
namely (\ref{genconda}) and (\ref{gencondb}) above respectively .
Hence, by Theorem \ref{cor3},
the map
$$\Psi_0:{\complex}^2\rightarrow {\complex}^2,
(z_1, z_2)\mapsto \left((1+2mV)^{\frac{1}{2}}e^{m(V-U)}z_1,
(1+2mU)^{\frac{1}{2}}e^{m(U-V)} z_2   \right)$$
is a special global symplectomorphism from  $({\C}^2, \omega_{m})$ into
$({\C}^2, \omega_{0})$.
\end{Example}

\begin{remar}\label{remarbefore} \rm
Notice that for $m>0$
 we cannot apply McDuff's  Theorem F in the Introduction in order to
 get the existence of
global symplectic coordinates on $({\complex}^2, \omega_m)$.
Indeed,  the sectional curvature of $({\C}^2,
g_{m})$ (where $g_m$ is the \K\ metric associated to $\omega_m$)
is positive at some point since $g_m$
is Ricci-flat but not flat.
\end{remar}

\begin{remar}\label{remarmcduff} \rm
In a forthcoming paper, among other properties,  we prove that
$({\complex}^2, g_m)$ cannot admit a \K\ immersion into any
complex space form, for all $m>\frac{1}{2}$
(this is achieved by using Calabi's diastasis function).
Therefore, the previous example
shows that the assumption in Theorem \ref{teor2} that $(M,
\omega_{\Phi})$ admits a \K\ immersion into some complex space
form is a sufficient but not a  necessary condition for the
existence of a special symplectic immersion into $({\complex}^n ,
\omega_0)$ .

\end{remar}

\small{}

\end{document}